\documentclass[11pt]{article}
\usepackage[french]{babel}
\usepackage[ansinew]{inputenc}
\usepackage{enumerate}
\usepackage{amsfonts}
\usepackage{amssymb}
\usepackage{amsmath}
\usepackage{amscd}
\usepackage{xypic}
\newtheorem{theo}{Th\'eor\`eme}[section]

\newtheorem{lem}[theo]{Lemme}
\newtheorem{cor}[theo]{Corollaire}
\newtheorem{rem}[theo]{Remarque}
\newtheorem{ex}[theo]{Exemple}

\newtheorem{propd}[theo]{Proposition et D\'efinition}
\def \endproof{\hfill$\square$ \bigskip}
\def\id{\textrm{\scriptsize id}}
\def\o{\textrm{O}}
\def\K{\mathbb{K}}
\def\D1{
\begin{equation*}
\raisebox{1.6cm}
{\xymatrix{ 
\ar[r]^-\id         & (\K,\o) \ar[r]^-\id \ar[d]^{\Gamma_i} & (\K,\o) \ar[r]^-\id \ar[d]^{\Gamma_{i-1}} & \cdots \ar[r]^-\id     & (\K,\o) \ar[r]^-\id \ar[d]^{\Gamma_1} & (\K,\o) \ar[d]^{\Gamma_0} \\
\ar[r]^-{\Pi_{i+1}} & (E_i,\o_i) \ar[r]^-{\Pi_i}            & (E_{i-1},\o_{i-1}) \ar[r]^-{\Pi_{i-1}}    & \cdots \ar[r]^-{\Pi_2} & (E_1,\o_1) \ar[r]^-{\Pi_1}            & (E_0,\o_0)                \\
\ar[r]^-{\omega_{i+1}}   & (Z_i,\o_i) \ar@{^{(}->}[u]_{J_i} \ar[r]^-{\omega_i} & (Z_{i-1},\o_{i-1}) \ar@{^{(}->}[u]_{J_{i-1}} \ar[r]^-{\omega_{i-1}}      & \cdots \ar[r]^-{\omega_2}   & (Z_1,\o_1) \ar@{^{(}->}[u]_{J_1} \ar[r]^-{\omega_1}              & (Z_0,\o_0) \ar@{^{(}->}[u]_{J_0}
}}
\end{equation*}
}

\begin{document}
\title{Une note \`a propos du Jacobien de $n$ fonctions holomorphes \`a l'origine de $\mathbb{C}^n$}
\author{M. Hickel}
\date{ }\maketitle
\renewcommand{\abstractname}{Abstract}
\begin{abstract}
 Let $f_{1},\ldots,f_{n}$ be $n$ germs of holomorphic functions at the origin of  $\mathbb{C}^{n}$, such that $f_{i}(0)=0$,   $1\leq i\leq n$ . We give a proof based on the  J. Lipman's theory of Residues via Hochschild Homology [L] that the jacobian of  $f_{1},\ldots,f_{n}$  belongs to the ideal generated by  $f_{1},\ldots,f_{n}$ if and only if the dimension of the germ of common zeros of  $f_{1},\ldots,f_{n}$ is strictly positive. In fact, we prove much more general results which are relatives versions of this result replacing the field $\mathbb{C}$ by convenients noetherian rings $\mathbf{A}$ (c.f Th. 3.1 and Th. 3.3). We then show a \L ojasiewicz inequality for the Jacobian analogous to the classical one by S. \L ojasiewicz for the gradient.
\footnote{1991 \it Mathematics Subject Classification; Primary 13A10; Secondary 13C10, 13J07.
 }\end{abstract}

 \section{Introduction}
  Soit $\mathcal {O}_{n}=\mathbb{C}\{ z_{1},\ldots,z_{n}\}$ l'anneau des germes de fonctions holomorphes \`a l'origine de 
$\mathbb{C}^{n}$ . Il est classique que si $f_{1},\ldots,f_{n}$  sont $n$ \'el\'ements de $\mathcal{O}_{n}$, ayant un z\'ero
isol\'e \`a l'origine de $\mathbb{C}^{n}$, alors leur jacobien $J(f_{1},\ldots,f_{n})=det(\partial  f_{i}/\partial  z_{j})$ n'appartient pas \`a l'id\'eal  $I$ de  $\mathcal{O}_{n}$  engendr\'e par  $f_{1},\ldots,f_{n}$. Motiv\'e par un r\'esultat de S. Spodzieza [S] prouvant un \'enonc\'e de E. Netto ([N] p.142) dans le cas de polyn\^omes homog\`enes, on peut se demander, comme le faisait A. Ploski [P], si la r\'eciproque de cette assertion est vraie.  Dans [H1] nous r\'epondions  par
l'affirmative \`a cette question.  Ce faisant, nous ne faisions que red\'ecouvrir dans ce contexte un r\'esultat du \`a W.V. Vasconcellos [W1] dans le cas polynomial. Nous pr\'esentons ici  dans un premier temps une preuve tr\`es simple de ce r\'esultat comme une cons\'equence de la th\'eorie des r\'esidus de J. Lipman [L], et qui vise \`a montrer l'efficacit\'e de celle-ci (Comparer avec les commentaires de nature historique de  [W2] p.268 \guillemotleft Scheja-Storch Formula \guillemotright et [W3] p. 129 \guillemotleft Tate Formula \guillemotright). En fait la th\'eorie des r\'esidus de J. Lipman via l'homologie de Hochschild nous permet de d\'emontrer des versions relatives de ces r\'esultats (c.f. Th. 3.1 et 3.3) et de remplacer le corps $\mathbb{C}$ par des anneaux noeth\'eriens $\mathbf{A}$ convenables. Si $(f_1,\ldots,f_n)$ ont un z\'ero isolé \`a l'origine, nous montrons dans un second temps que, d\'esignant par $m$ le rang en $0$ de la matrice Jacobienne $\frac{\partial f_i}{\partial z_j}(0)$ on a : si $n-m\geq 3$ $J(f_1,\ldots,f_n) \in \overline{I^{\theta }}$, o\`u $\theta $ est un nombre rationnel \underline{strictement} plus grand que $1$ et que si $n-m=2$ on a : $J(f_1,\ldots,f_n) \in \overline{I}$ o\`u $\overline{I^q}$ d\'esigne la cl\^oture int\'egrale des id\'eaux (c.f. plus bas et [H-S], [L-T] pour des rappels concernant ces notions). On am\'eliore ainsi un r\'esultat de B. Teissier [T1] et on prouve  une in\'egalit\'e de \L ojasiewicz pour le Jacobien, analogue \`a celle du Gradient  [Lo], [B-M],[T2]. Tout d'abord l'anneau $\mathcal{O}_n$ \'etant local noeth\'erien, son compl\'et\'e $\hat{\mathcal{O}}_n=\mathbb{C}[[X_1,\ldots,X_n]]$ est fid\`element plat sur $\mathcal{O}_n$ et donc pour tout id\'eal $I\subset \mathcal{O}_n$, on a : $I.\hat{\mathcal{O}}_n\cap \mathcal{O}_n=I$. Le passage au compl\'et\'e ne changeant pas les dimensions, il suffira de travailler dans $\mathbb{C}[[X_1,\ldots,X_n]]$. On a alors : 
 
\begin{theo}\textrm{ }\\
Soient $k$ un corps de caract\'eristique z\'ero et $k[[X_{1},\ldots,X_{n}]]$ l'anneau des s\'eries formelles en $n$
variables \`a coefficients dans $k$. Pour $f_{1},\ldots,f_{n}$ $\in  \ k[[X_{1},\ldots,X_{n}]]$
d\'efinissant un id\'eal propre, les propri\'et\'es suivantes sont \'equivalentes :
\begin{itemize}
\item[1)] $Dim\ k{\rm [[}X{\rm ]]}/(f_{1},\ldots,f_{n})\ >\ 0$
\rm(\it o\`u Dim d\'esigne la dimension de Krull de l'anneau consid\'er\'e).
\item[2)]Le jacobien $J_{X}(f_{1},\ldots,f_{n})$ est dans l'id\'eal $(f_{1},\ldots,f_{n})k{\rm [[}X{\rm
]]}$.
\end{itemize}
De plus si $dim_k \frac{k[[X_1,\ldots,X_n]]}{(f_1,\ldots,f_n)}<+\infty$ (comme $k$-espace vectoriel) i.e. $Dim \frac{k[[X_1,\ldots,X_n]]}{(f_1,\ldots,f_n)}$ $=0$, on a : $$J_X(f_1,\ldots,f_n).(X_1,\ldots,X_n)\subset (f_1,\ldots,f_n).$$ 
Autrement dit la classe de $J_X(f_1,\ldots,f_n)$ dans $\frac{k[[X_1,\ldots,X_n]]}{(f_1,\ldots,f_n)}$ est un g\'en\'erateur du socle de cette alg\`ebre Artinienne.
\end{theo} 
La derni\`ere assertion est aussi valide lorsque $k$ est de caract\'eristique positive $p$, et $p$ ne divise pas $dim_k \frac{k[[X_1,\ldots,X_n]]}{(f_1,\ldots,f_n)}$. Cette assertion est due \`a S. Scheja, U. Storch [S-S]. L'implication $1) \Longrightarrow 2)$ est due \`a W.V. Vasconcelos dans le cas des alg\`ebres de polyn\^omes [W1]. Nous pr\'esentons une preuve  d'une version relative de toutes ces assertions via la th\'eorie des r\'esidus de J. Lipman [L] (c.f Th.3.1, Th.3.3 et Cor. 3.5) rempla\c{c}ant le corps $k$ par des anneaux noeth\'eriens convenables.  Nous rappelons \`a la section 2 quelques \'el\'ements concernant cette th\'eorie et prouvons les r\'esultats \`a la section 3). \\
 Nous prouvons ensuite le r\'esultat suivant :
\begin{theo}\textrm{ }\\
Soient $k$ un corps  de caract\'eristique 0 et $f_1,\ldots,f_n \in k[[X_1,\ldots,X_n]]$ tels que $f_i(0)=0$ et l'id\'eal $I=(f_1,\ldots,f_n)$ est $\mathfrak{m}=(X_1,\ldots,X_n)$ primaire (Ainsi $J_X(f)\not \in I$). D\'esignons par $m$ le rang en $0$ de la matrice Jacobienne $\frac{\partial f_i}{\partial X_j}(0)$ :
\begin{itemize}
\item[1)] Si $n-m=2$, on a $J_X(f_1,\ldots,f_n) \in \overline{I}$.
\item[2)] Si $n-m\geq 3$, il existe un nombre rationnel $\theta=p/q >1$ tel que $J_X(f_1,\ldots,f_n) \in \overline{I^{\theta} }$ o\`u encore $J_X(f_1,\ldots,f_n)^q \in \overline{I^p}$.
 \item[3)] Si $k$ est muni d'une valeur absolue non triviale $\vert \quad \vert$: $k\longrightarrow [0,+\infty[$ et si les $f_i$ sont des germes de s\'eries convergentes alors :
\item [-] Dans le cas 1), $\exists r>0,\hspace{3pt}\exists C>0 / \hspace{3pt}\parallel x\parallel \leq r\Longrightarrow \vert J_X(f_1,\ldots,f_n)(x)\vert \leq C (\sum_{i=1}^n\vert f_i(x)\vert).$
\item [-] Dans le cas 2), $\exists r>0,\hspace{3pt}\exists C>0 / \hspace{3pt}\parallel x\parallel \leq r\Longrightarrow \vert J_X(f_1,\ldots,f_n)(x)\vert \leq C (\sum_{i=1}^n\vert f_i(x)\vert)^{\theta }$, o\`u $\theta $ est un rationnel strictement plus grand que 1.
\end{itemize}
\end{theo}
Bien entendu, il n'y a pas de r\'esultats similaires lorsque $n-m=1$, car on est alors ramen\'e au cas o\`u $n=1$. La conclusion de 2) n'est pas valide en g\'en\'eral dans le cas $n-m=2$ comme nous le verrons par des exemples. Nous rappelons quelques \'el\'ements concernant la cl\^oture int\'egrale des id\'eaux \`a la section 4) et y prouvons le r\'esultat ci-dessus.

\section{Quelques rappels sur la th\'eorie des r\'esidus de J. Lipman [L]}
On rappelle ici quelques d\'efinitions et r\'esultats de [L] § 3, pour d'autres pr\'ecisions et applications nous renvoyons le lecteur \`a [L], [BO], [B-H1], [B-H2]. On consid\`ere $\mathbf{R}$ une $\mathbf{A}$-alg\`ebre commutative et unitaire de morphisme structural $h$ : $\mathbf{A}\longrightarrow \mathbf{R}$. On consid\`ere $f=(f_1,\ldots,f_n)$ une suite quasi-r\'eguli\`ere de $\mathbf{R}$, $I$ l'id\'eal qu'elle engendre, et $\hat{\mathbf{R}}$ le compl\'et\'e $I$-adique de $\mathbf{R}$. Rappelons que $f$ est dite quasi-r\'eguli\`ere si $f.\mathbf{R}=I.\mathbf{R}\neq \mathbf{R}$ et si le $\mathbf{R}/I$ homomorphisme de $(\mathbf{R}/I)[X_1,\ldots,X_n]\longrightarrow gr_{I}(\mathbf{R})=\bigoplus_{m\geq 0}\frac{I^m}{I^{m+1}}$ qui \`a un polyn\^ome homog\`ene $P_m(X_1,\ldots,X_n)$ de degr\'e $m$ fait correspondre la classe de $P_m(f_1,\ldots,f_n)$ dans $\frac{I^m}{I^{m+1}}$ est un isomorphisme. On suppose que $\mathbf{P}=\mathbf{R}/I$ est un $\mathbf{A}$-module projectif de type fini et l'on note $\sigma :$ $\mathbf{R}/I\longrightarrow \mathbf{R}$ une section $\mathbf{A}$ lin\'eaire de la projection canonique de $\mathbf{R}\longrightarrow \frac{\mathbf{R}}{I}$. Tout \'el\'ement $r$ de $\mathbf{R}$ admet dans $\hat{\mathbf{R}}$ une \'ecriture unique de la forme :
$$r=\sum_{\alpha \in \mathbb{N}^n}\sigma (r_{\alpha }).f^{\alpha }, \textit{ avec } r_{\alpha }\in \mathbf{P}$$
Par cons\'equent, il existe des uniques $\gamma _{\alpha }\in Hom_{\mathbf{A}}(\mathbf{P},\mathbf{P})$ tels que :
$$\forall p \in \mathbf{P},\hspace{3pt}r.\sigma (p)=\sum_{\alpha \in \mathbb{N}^n}\sigma (\gamma _{\alpha }(p)).f^{\alpha }$$
On d\'efinit l'\'el\'ement $r^{\sharp}\in Hom_{\mathbf{A}}(\mathbf{P},\mathbf{P})[[f]]$ par :
$$ r^{\sharp}=\sum_{\alpha \in \mathbb{N}^n}\gamma _\alpha f^{\alpha }$$
Soit $\epsilon _i$ le n-uplet dont toutes les composantes sont nulles sauf la $i^{\grave{e}me}$ qui vaut $1$. On consid\`ere l'\'el\'ement suivant de $Hom_{\mathbf{A}}(\mathbf{P},\mathbf{P})[[f]]$ :
$$\frac{\partial }{\partial f_i}(r^{\sharp})=\sum_{\alpha =(\alpha _1,\ldots,\alpha _n)\in \mathbb{N}^n}\alpha _i\gamma _{\alpha}.f^{\alpha -\epsilon _i}$$
$\mathbf{P}$ \'etant une $\mathbf{A}$-alg\`ebre projective de type fini, on dispose d'un op\'erateur trace canonique $Tr_{\mathbf{P}/\mathbf{A}} \in hom_{\mathbf{A}}(End_{\mathbf{A}}(\mathbf{P}),\mathbf{A})$ qui \`a tout \'el\'ement de $End_{\mathbf{A}}(\mathbf{P})$ fait correspondre sa trace (c.f. [B1] chap. II §4 p. 78-79). On a alors :
\begin{theo} ([L] §3 corollaire (3.7))\\
Soient $f=(f_1,\ldots,f_n)$ une suite quasi-r\'eguli\`ere de $\mathbf{R}$ telle que le $\mathbf{P}=\mathbf{R}/I$ module soit un  $\mathbf{A}$-module projectif de type fini, $\omega =rdr_1\wedge \ldots \wedge dr_n \in \wedge ^q \Omega _{\mathbf{R}/\mathbf{A}}$, o\`u  $\Omega _{\mathbf{R}/\mathbf{A}}$ est le $\mathbf{R}$-module des diff\'erentielles de la $\mathbf{A}$-alg\`ebre $\mathbf{R}$.  Pour tout multi-indice $(m_1,\ldots,m_n)\in (\mathbb{N}^{*})^n$, posant :
$$r^{\sharp}det(\frac{\partial }{\partial f_i}(r_j^{\sharp}))=\sum_{\alpha \in \mathbb{N}^n}\delta _\alpha .f^\alpha \in hom_{\mathbf{A}}(\mathbf{P},\mathbf{P})[[f]]$$
On a :
$$Res_{\mathbf{R}/\mathbf{A}}\begin{bmatrix} & rdr_1\wedge \ldots \wedge dr_n & \\
& f_1^{m_1},\ldots,f_n^{m_n}& \end{bmatrix}=Tr_{\mathbf{P}/\mathbf{A}}(\delta _{m_1-1,\ldots,m_n-1})$$
\end{theo} 

Ce th\'eor\`eme peut servir de d\'efinition, dans le cas consid\'er\'e, pour le symbole r\'esidu $Res_{\mathbf{R}/\mathbf{A}}$ de J. Lipman. 
\begin{cor} ([L] corollaire 3.8 p.51)\\
Soient $\mathbf{A}$, $\mathbf{R}$ ,$f=(f_1,\ldots,f_n)$, $I$ et $\mathbf{P}=\mathbf{R}/I$ comme ci-dessus alors pour tout $r \in \mathbf{R}$ :
$$Res_{\mathbf{R}/\mathbf{A}}\begin{bmatrix} & rdf_1\wedge \ldots \wedge df_n & \\
& f_1,\ldots,f_n& \end{bmatrix}=Tr_{\mathbf{P}/\mathbf{A}}(p_r)$$
o\`u $p_r$ est l'op\'erateur de multiplication par $r$ dans $\mathbf{P}\longrightarrow \mathbf{P}$,  $p\longrightarrow \overline{r}.p$
\end{cor}

$\mathbf{A}$ \'etant un anneau commutatif unitaire noeth\'erien, consid\'erons la situation particuli\`ere o\`u $\mathbf{R}=\mathbf{A}[X_1,\ldots,X_n]$ ou bien $\mathbf{A}[[X_1,\ldots,X_n]]$ ou encore le localis\'e en un id\'eal premier d'un tel anneau. D\'esignons toujours par $f$ une suite quasi-r\'eguli\`ere telle que le $\mathbf{P}=\mathbf{R}/I$ soit un $\mathbf{A}$-module projectif de type fini. On a alors :
\begin{cor}(c.f. [BO] Chap.1 prop. 2.2.2 p. 21, prop. 2.3.1 p.29 , prop. 4.1.2 p. 48)\\
$Hom_{\mathbf{A}}(\mathbf{P},\mathbf{A})$ est un $\mathbf{P}$-module libre de rang $1$ dont un g\'en\'erateur est l'homomorphisme r\'esidu : $$r \longrightarrow Res_{\mathbf{R}/\mathbf{A}}\begin{bmatrix} & rdX_1\wedge \ldots \wedge dX_n & \\
& f_1,\ldots,f_n& \end{bmatrix}.$$ 
En particulier, on a la propri\'et\'e de non d\'eg\'en\'erescence suivante :
$$r \in I=(f_1,\ldots,f_n)\Longleftrightarrow \forall h \in \mathbf{R},\hspace{3pt} Res_{\mathbf{R}/\mathbf{A}}\begin{bmatrix} & rh dX_1\wedge \ldots \wedge dX_n & \\
& f_1,\ldots,f_n& \end{bmatrix}=0$$
\end{cor}
 \section{Versions relatives  de 1.1}
Nous pouvons \`a pr\'esent \'enoncer le r\'esultat suivant qui est une version relative de $2)\Longrightarrow 1)$ et de la derni\`ere assertion de 1.1.
\begin{theo}\textrm{ }\\
Soit $\mathbf{A}$ un anneau commutatif noeth\'erien int\`egre de caract\'eristique $p\geq 0$. D\'esignons par $\mathbf{R}$ l'anneau $\mathbf{A}[X_1,\ldots,X_n]$ ou $\mathbf{A}[[X_1,\ldots,X_n]]$ ou bien le localis\'e en un id\'eal premier d'un de ces deux anneaux. Soit $f=(f_1,\ldots,f_n)\in \mathbf{R}$ une suite d'\'el\'ements de $\mathbf{R}$ telle que $I=(f_1,\ldots,f_n).\mathbf{R}\neq \mathbf{R}$. On suppose que :\\
- a) $f$ est quasi-r\'eguli\`ere\\
- b) $\mathbf{P}=\mathbf{R}/I$ est un $\mathbf{A}$ module projectif de type fini dont on note $r$ le rang.\\
Alors notant par $J_X(f_1,\ldots,f_n)=det(\frac{\partial f_i}{\partial X_j})$ le Jacobien de $(f_1,\ldots,f_n)$ on a :\\
1) $J_X(f_1,\ldots,f_n).\sqrt{I}\subset I$\\
2) Si $r$ est premier \`a $p$ (en particulier si $p=0$) alors $J_X(f_1,\ldots,f_n) \not \in I$.\\
 Les conclusions ci-dessus 1) et 2) sont aussi satisfaites sous les hypoth\`eses alternatives suivantes. $\mathbf{A}$ est un anneau local noeth\'erien d'\'egale caract\'eristique $p$, de dimension $0$ (non n\'ecessairement int\`egre), et $a)$ est satisfaite et $b)$ est remplac\'ee par :\\
- b') $\mathbf{P}=\mathbf{R}/I$ est un $\mathbf{A}$ module libre de type fini  et la longueur du $\mathbf{A}$-module $\mathbf{P}$ (qui est n\'ecessairement finie) est premi\`ere \`a $p$.
\end{theo}
Avant de prouver le th\'eor\`eme notons que dans le cas particulier o\`u $A$ est un corps $k$ infini les hypoth\`eses $a)$ et $b)$ se r\'eduisent \`a $dim_{k}k[X_1,\ldots,X_n]/I<+\infty$ ou $dim_k k[[X_1,\ldots,X_n]]/I<+\infty$. En effet cette seule condition garantit la quasi-r\'egularit\'e de la suite $f=(f_1,\ldots,f_n)$ (c.f. [Bo] appendice prop. 3.1 et cor. 1.9). Ceci \'equivaut au fait que $Dim k[X_1,\ldots,X_n]/I=0$ ou $Dim k[[X_1,\ldots,X_n]]/I=0$ et l'on retrouve les r\'esultats de 1.1 annonc\'es.\\

\noindent \textit{Preuve :}\\
Soit $\sigma \in End_{\mathbf{A}}(\mathbf{P})$, $\sigma :$ $\mathbf{P}\longrightarrow \mathbf{P}$. Puisque $\mathbf{P}$ est un $\mathbf{A}$-module projectif de type fini, on peut trouver un entier $n$ et un morphisme surjectif $u$ : $\mathbf{A}^n\longrightarrow \mathbf{P}\longrightarrow 0$. Puisque $\mathbf{P}$ est projectif, on peut trouver un morphisme $v$ : $\mathbf{P}\longrightarrow \mathbf{A}^n$ tel que $v$ soit une section de $u$, c'est \`a dire :
$$ u\circ v=id_{\mathbf{P}}.$$
Consid\'erons l'endomorphisme de $\mathbf{A}^n$, $l_{\sigma }\in End_{\mathbf{A}}(\mathbf{A}^n)$ d\'efini par :
$$ l_{\sigma }=v\circ \sigma \circ u.$$
On a alors le lemme suivant:
\begin{lem}\textrm{ }\\
Soient $\sigma \in End_{\mathbf{A}}(\mathbf{P})$ et $l_{\sigma }\in End_{\mathbf{A}}(\mathbf{A}^n)$ comme ci-dessus.\\
1) Si $\sigma $ est nilpotent, $l_{\sigma }$ est nilpotent.\\
2) $Tr_{\mathbf{P}/\mathbf{A}}(\sigma )=Tr_{\mathbf{A}^n/\mathbf{A}}(l_{\sigma })$\\
3) Si $\sigma $ est nilpotent, $Tr_{\mathbf{P}/\mathbf{A}}(\sigma )=0$.
\end{lem}
\noindent \textit{Preuve :}\\
1) Supposons $\sigma ^m=0$ pour un entier $m$. Alors $l_{\sigma }^m=(v \circ \sigma \circ u)\circ \ldots \circ (v \circ \sigma \circ u)$. Comme $u\circ v=id_{\mathbf{P}}$, on a $l_{\sigma }^m=v\circ \sigma ^m\circ u$. D'o\`u $l_{\sigma} ^m=0$,  puisque $\sigma ^m=0$.\\
2) Notons $\mathbf{P}^{*}=Hom_{\mathbf{A}}(\mathbf{P},\mathbf{A})$. D'apr\`es [B1] et l'isomorphisme canonique :
$$\mathbf{P}^*\otimes _{\mathbf{A}}\mathbf{P}\longrightarrow End_{\mathbf{A}}(\mathbf{P})$$
on peut trouver des $x_i^{*} \in \mathbf{P}^{*}$, des $y_i \in \mathbf{P}$ tels que pour tout $x \in \mathbf{P}$ :
$$\sigma (x)=\sum_{i=1}^m <x,x_i^{*}>.y_i$$
On a ainsi $Tr_{\mathbf{P}/\mathbf{A}}(\sigma )=\sum_{i=1}^m <y_i,x_{i}^{*}>$ ([B1] Chap II § 4 78-79). Soit $(e_i)_{1\leq i \leq n}$ la base canonique de $\mathbf{A}^n$. Posons :
$$a_i=u(e_i),\hspace{3pt}\sigma (a_i)=\sum_{j=1}^m<a_i,x_{j}^{*}>.y_j \textrm{ et }v(y_j)=\sum_{k=1}^nv_{j,k}.e_k$$
Ainsi :
\begin{align}
l_{\sigma }(e_i) & =v(\sigma (a_i))=\sum_{j=1}^m <a_i,x_{j}^{*}>.v(y_j)\\
      & =\sum_{j=1}^m<a_i,x_j^{*}>(\sum_{k=1}^nv_{j,k}.e_k)\\
      & =\sum_{k=1}^n(\sum_{j=1}^m<a_i,x_j^{*}>v_{j,k}).e_k
\end{align}
Ainsi :
\begin{align}
Tr_{\mathbf{A}^n/\mathbf{A}}(l_{\sigma }) & =\sum_{i=1}^n (\sum_{j=1}^m<a_i,x_j^{*}>.v_{j,i})\\
  & = \sum_{j=1}^m<\sum_{i=1}^nv_{j,i}a_i,x_j^{*}>\\
 & = \sum_{j=1}^m<\sum_{i=1}^nv_{j,i}u(e_i),x_j^{*}>\\
 & =\sum_{j=1}^m<u(\sum_{i=1}^nv_{j,i}.e_i),x_j^{*}>\\
 & = \sum_{j=1}^m<u\circ v(y_j),x_j^*>=\sum_{j=1}^m<y_j,x_j^*>=Tr_{\mathbf{P}/\mathbf{A}}(\sigma )
\end{align}
3) Soient $K=frac(A)$ et $l'_{\sigma }$ : $K^n \longrightarrow K^n$ obtenu \`a partir de $l_{\sigma }$ par extension des scalaires. Notons $\overline{K}$ une cl\^oture alg\'ebrique de $K$ et $\overline{l_{\sigma }}$ l'extension \`a $\overline{K}^n$ de $l'_{\sigma }$. Comme $\sigma $ est nilpotent par 1), $l_{\sigma }$ est nilpotent et par cons\'equent $l'_{\sigma }$ et $\overline{l_{\sigma }}$ sont nilpotents. Maintenant $\overline{l_{\sigma }}$ : $\overline{K}^n\longrightarrow \overline{K}^n$ \'etant une application lin\'eaire nilpotente toutes ses valeurs propres sont nulles et l'existence de la r\'eduction des matrices carr\'ees \`a coefficients dans $\overline{K}$ nous assure que $Tr_{\overline{K}^n/\overline{K}}(\overline{l_{\sigma }})=0$. Mais d'apr\`es 2) $Tr_{\mathbf{P}/\mathbf{A}}(\sigma )=Tr_{\mathbf{A}^n/\mathbf{A}}(l_{\sigma })=Tr_{\overline{K}^n/\overline{K}}(\overline{l_{\sigma }})=0$. \endproof

Nous pouvons \`a pr\'esent prouver 3.1. Concernant la premi\`ere affirmation, soit $r \in \sqrt{I}$. Il s'agit de voir que $J_X(f).r \in I$. Pour cela, il est \'equivalent par 2.3 de prouver que pour tout $h\in \mathbf{R}$, on a :
 $$Res_{\mathbf{R}/\mathbf{A}}\begin{bmatrix} & rJ_X(f)h dX_1\wedge \ldots \wedge dX_n & \\
& f_1,\ldots,f_n& \end{bmatrix}=0$$
Or 
$$ Res_{\mathbf{R}/\mathbf{A}}\begin{bmatrix} & rJ_X(f)h dX_1\wedge \ldots \wedge dX_n & \\
& f_1,\ldots,f_n& \end{bmatrix}= Res_{\mathbf{R}/\mathbf{A}}\begin{bmatrix} & rh df_1\wedge \ldots \wedge df_n & \\
& f_1,\ldots,f_n& \end{bmatrix}$$
Donc en vertu de 2.2 ce dernier r\'esidu n'est autre que la trace de l'op\'erateur $p_{hr}$ de multiplication par $hr$ dans $\mathbf{P}$. Comme $r\in \sqrt{I}$, $rh\in \sqrt{I}$ et donc $p_{rh}$ est un op\'erateur nilpotent. Le lemme $3.2$ nous assure que sa trace est nulle et nous permet donc de conclure \`a la premi\`ere affirmation de 3.1.\\
Concernant la seconde assertion de $3.1$. Supposons que $J_X(f_1,\ldots,f_n)\in I$. Alors pour tout id\'eal maximal $\mathfrak{m}$ de $\mathbf{A}$, d\'esignant toujours par $J_X(f)$ la classe de $J_X(f)$ dans le localis\'e $\mathbf{R}_{\mathfrak{m}}$ de $\mathbf{R}$ en $\mathfrak{m}$, on aurait $J_X(f)\in I_{\mathfrak{m}}$. On aurait ainsi :
$$Res_{\mathbf{R}_{\mathfrak{m}}/\mathbf{A}_{\mathfrak{m}}}\begin{bmatrix} & J_X(f) dX_1\wedge \ldots \wedge dX_n & \\
& f_1,\ldots,f_n& \end{bmatrix}= Res_{\mathbf{R}_{\mathfrak{m}}/\mathbf{A}_{\mathfrak{m}}}\begin{bmatrix} &  df_1\wedge \ldots \wedge df_n & \\
& f_1,\ldots,f_n& \end{bmatrix}=0$$
D'apr\`es 2.2 ce dernier r\'esidu n'est autre que la trace de l'op\'erateur identit\'e de $\mathbf{P}_{\mathfrak{m}}$. Comme $\mathbf{P}_{\mathfrak{m}}$ est un $\mathbf{A}_{\mathfrak{m}}$ module libre de rang $r$ (c.f [B2] II § 5 p.141), ce r\'esidu vaut donc $r$. Il ne  peut donc \^etre nul que si la caract\'eristique de $\mathbf{A}_{\mathfrak{m}}$, c'est \`a dire celle de $\mathbf{A}$, divise $r$. \\
Nous indiquons maintenant comment modifier la preuve ci-dessus dans les hypoth\`eses alternatives. La seule chose qui soit modifi\'ee est le fait que l'on utilise des arguments diff\'erents pour faire les calculs de trace. Soient $\hat{\mathbf{A}}$ et $\hat{\mathbf{P}}$ les compl\'et\'es $\mathfrak{m}_{\mathbf{A}}$-adiques de $\mathbf{A}$ et $\mathbf{P}$. Comme $\mathbf{P}$ est suppos\'e \^etre un $\mathbf{A}$-module libre, $\hat{\mathbf{P}}$ est un $\hat{\mathbf{A}}$ libre (car le morphisme $\mathbf{A}\longrightarrow \hat{\mathbf{A}}$ est fid\`element plat). Si $g$ est un \'el\'ement de $\mathbf{P}$ la matrice de l'op\'erateur de $p_g$ de multiplication par $g$ comme \'el\'ement de $End_{\mathbf{A}}(\mathbf{P})$ est la m\^eme que celle de l'op\'erareur de multiplication $\hat{p_g}$ comme \'el\'ement de $End_{\hat{\mathbf{A}}}(\hat{\mathbf{P}})$. Il suffit pour le voir de consid\'erer un isomorphisme d'un $\mathbf{A}^m$ sur $\mathbf{P}$ et de tensoriser par $\hat{\mathbf{A}}$. Ainsi :
$$ Tr_{\mathbf{P}/\mathbf{A}}(p_g)=Tr_{\hat{\mathbf{P}}/\hat{\mathbf{A}}}(\hat{p_g})$$
Ensuite notons $K=\frac{\mathbf{A}}{\mathfrak{m}_{\mathbf{A}}}$, le th\'eor\`eme de structure des anneaux locaux complets d'\'egale caract\'eristique de I. Cohen nous dit que $\hat{\mathbf{A}}$ est isomorphe \`a $K[[Y]]/J$, o\`u $J$ est un id\'eal   $(Y)$ primaire puisque $\mathbf{A}$ et donc $\hat{\mathbf{A}}$ sont de dimension $0$. Ainsi $\hat{\mathbf{P}}$ est une $K$-alg\`ebre et un $K$-espace vectoriel de dimension finie, dimension qui n'est autre que la longueur du $A$-module Artinien $\mathbf{P}$. Maintenant on a (c.f. [B1] chap III § 9 p.110 (21)):
$$Tr_{\hat{\mathbf{P}}/\hat{\mathbf{A}}}(\hat{p_g})=Tr_{\hat{\mathbf{P}}/K}(\hat{p_g})$$
Dans le cas o\`u $p_g$ est nilpotente, cette trace est donc nulle puisque $\hat{p_g}$ est un endomorphisme nilpotent sur un $K$-espace vectoriel de dimension finie. Et si $g=1$ cette trace est la classe dans $K$, de la longueur qui intervient dans l'\`enonc\'e ou encore celle de $dim_K\hat{\mathbf{P}}$. Cette classe est non nulle puisque par hypoth\`ese ce nombre est premier \`a la caract\'eristique de $K$. Ceci permet de conclure comme plus haut aux assertions 1) et 2) de 3.1).
\endproof

Nous allons maintenant \'enoncer une version relative de $1)\Longrightarrow 2)$ dans 1.1. Pour la notion d'anneau de Gorenstein nous renvoyons \`a [E] Chap. 21 p. 525 ou [B3] § 3, p. 48.
 
\begin{theo}\textrm{ }\\
Soit $\mathbf{A}$ un anneau local de Gorenstein, Artinien d'\'egale caract\'eristique z\'ero. Posons $\mathbf{R}=\mathbf{A}[[X_1,\ldots,X_n]]$ et $(f_1,\ldots,f_n)\in \mathbf{R}$ tels que $(f_1,\ldots,f_n).\mathbf{R}=I. \mathbf{R}\neq \mathbf{R}$. Alors si $\mathbf{R}/I$ n'est pas Artinien alors :
$$J_X(f_1,\ldots,f_n)\in I$$
\end{theo}
\noindent \textit{Preuve :}\\
Remarquons d'abord que $\mathbf{A}$ \'etant de Gorenstein, il est de Cohen-Macauley. Il en est de m\^eme de $\mathbf{R}$ qui est de Gorenstein et de Cohen-Macauley. Comme la dimension de Krull de $\mathbf{R}$ est $n$, il s'en suit que toute suite $(g_1,\ldots,g_n)$ qui engendre un id\'eal $\mathfrak{m}_{\mathbf{R}}$ primaire est r\'eguli\`ere donc quasi-r\'eguli\`ere puisque l'on travaille ici dans un cadre local. Par ailleurs $\mathbf{R}/(g_1,\ldots,g_n).\mathbf{R}$ est un $\mathbf{A}$ module libre de type fini. En effet, $\mathbf{A}$ \'etant local Artinien il est s\'epar\'e et complet pour la topologie $\mathfrak{m}_{\mathbf{A}}$-adique. Par le th\'eor\`eme de structure de I. Cohen des anneaux locaux complets d'\'egale caract\'eristique on a $\mathbf{A}\simeq k[[Y]]/\mathfrak{a}$, o\`u $\mathfrak{a}$ est $\mathfrak{m}_{\mathbf{A}}$-primaire, $k\simeq \mathbf{A}/\mathfrak{m}_{\mathbf{A}}$ et $Y=(Y_1,\ldots,Y_p)$. Notons $g_1(Y,X),\ldots,g_n(Y,X)$ des \'el\'ements de $k[[Y,X]]$ repr\'esentant les $g_i \in \mathbf{R}\simeq k[[Y,X]]/\mathfrak{a}$. Puisque $\mathbf{R}/(\mathfrak{m}_{\mathbf{A}},g_1,\ldots,g_n).\mathbf{R}$ est Artinien, on peut affirmer que :
$$\frac{k[[X]]}{(g_1(0,X),\ldots,g_n(0,X))}\simeq \frac{k[[Y,X]]}{(Y,g_1(Y,X),\ldots,g_n(Y,X))}\simeq \frac{\mathbf{R}}{(\mathfrak{m}_{\mathbf{A}},g_1,\ldots,g_n)}$$
est un $k$-espace vectoriel de dimension finie. Il en d\'ecoule que $g_1(0,X),\ldots,g_n(0,X)$ est une suite r\'eguli\`ere de $k[[X]]\simeq \mathbf{R}/\mathfrak{m}_{\mathbf{A}}.\mathbf{R}$. Par cons\'equent, puisque $\mathbf{R}$ est $\mathbf{A}$-plat, le crit\`ere local de platitude (c.f. [M] Th. 22.5 et Cor. p. 176-177) nous assure que $\mathbf{R}/(g_1,\ldots,g_n)\mathbf{R}$ est un $\mathbf{A}$-module (de type fini) plat. Puisque $\mathbf{A}$ est local, $\mathbf{R}/(g_1,\ldots,g_n)\mathbf{R}$ est un $\mathbf{A}$-module libre de type fini. On pourra donc appliquer 3.1 \`a de telles suites $(g_1,\ldots,g_n)$. \\
Soit $(f)=(f_1,\ldots,f_n)$ satisfaisant les hypoth\`eses du th\'eor\`eme. Notons d'abord que le r\'esultat est trivial si le nombre minimal de g\'en\'erateurs de $I$ est strictement plus petit que $n$. En effet, dans cette \'eventualit\'e, si $u_1,\ldots,u_k$, $k<n$ est un syst\`eme de g\'en\'erateurs de $I$. On a :
$$f_i=\sum_{j=1}^na_{i,j}u_j,\hspace{3pt}d_Xf_i=\sum_{j=1}^na_{i,j}d_Xu_j \pmod{I},\hspace{3pt}a_{i,j}\in \mathbf{R}$$
Ainsi $J_X(f)=0 \pmod{I}$. On peut donc supposer sans restriction que le nombre minimal de g\'en\'erateurs de $I$ est $n$.
Notons $\mathfrak{m}_{\mathbf{R}}$ l'id\'eal maximal de $\mathbf{R}$. Par le th\'eor\`eme d'intersection de Krull, on a :
$$ I=\cap _{a \in \mathbb{N}}(I+\mathfrak{m}_{\mathbf{R}}^a)$$
Il suffit donc de d\'emontrer que : $J_X(f) \in I+\mathfrak{m}_{\mathbf{R}}^a$ pour $a$ assez grand. Fixons $a$ et soient     $n_a$ le nombre minimal de g\'en\'erateurs de $I_a$ et $g_1,\ldots,g_{n_a}$ un syst\`eme minimal de g\'en\'erateurs de $I_a$. Comme $Dim(\mathbf{R})=n$, les r\'esultats de Northcott-Rees [N-R] sur la r\'eduction des id\'eaux nous assurent que quitte \`a remplacer les $g_i$ par $n_a$ combinaisons lin\'eaires \`a coefficients dans $\mathbb{Q}$ suffisamment g\'en\'eraux des $g_i$, on peut supposer que toute sous-famille $g_{i_1},\ldots,g_{i_n}$ est telle que l'id\'eal $(g_{i_1},\ldots,g_{i_n})$ est $\mathfrak{m}_{\mathbf{R}}$-primaire. Deux \'eventualit\'es peuvent se pr\'esenter.\\
- Soit $n_a>n$. Alors :
$$f_i=\sum_{j=1}^{n_a}a_{i,j}.g_j,\textrm{ donc }d_Xf_i=\sum_{j=1}^{n_a}a_{i,j}.d_Xg_j \pmod{I_a}$$
Ainsi :
$$J_X(f)=\sum_{1\leq i_1<\ldots<i_n}A_{(i_1,\ldots,i_n)}J_X(g_{i_1},\ldots,g_{i_n}) \pmod{I_a}, \hspace{3pt}A_{(i_1,\ldots,i_n)}\in \mathbf{R}$$
Il suffit donc de voir que chaque $J_X(g_{i_1},\ldots,g_{i_n}) $ est dans $I_a$. Soient $\underline{i}=(i_1,\ldots,i_n)$ un tel multi-indice et $\mathbf{P}_{\underline{i}}=\mathbf{R}/(g_{i_1},\ldots,g_{i_n})$. Notons $S_{\underline{i}}$ le socle de $\mathbf{P}_{\underline{i}}$ i.e. $S_{\underline{i}}=\{x \in \mathbf{P}_{\underline{i}}/\mathfrak{m}_{\underline{i}}.x=0\}$ o\`u $\mathfrak{m}_{\underline{i}}$ est l'id\'eal maximal de $\mathbf{P}_{\underline{i}}$. D'apr\`es 3.1 $\overline{J}_X(g_{i_1},\ldots,g_{i_n})\in S_{\underline{i}}$ (car $g_{i_1},\ldots,g_{i_n}$ est une suite r\'eguli\`ere de $\mathbf{R}$). Maintenant on a les inclusions strictes :
$$0\subset I_a.\mathbf{P}_{\underline{i}}\subset \mathbf{P}_{\underline{i}}$$
$\mathbf{P}_{\underline{i}}$ \'etant Artinien, tout sous-module non trivial de $\mathbf{P}_{\underline{i}}$ contient un sous-module simple. Un sous-module simple est par le lemme de Nakayama inclus dans le socle de $\mathbf{P}_{\underline{i}}$. Maintenant $\mathbf{R}$ \'etant de Gorenstein et $g_{i_1},\ldots,g_{i_n}$ une suite r\'eguli\`ere de $\mathbf{R}$, $\mathbf{P}_{\underline{i}}$ est de Gorenstein (c.f [B3] exemple 2 p. 48) et donc son socle est simple et c'est le seul sous-module simple de $\mathbf{P}_{\underline{i}}$. On a donc :
$$\overline{J}_X(g_{i_1},\ldots,g_{i_n})\in S_{\underline{i}}\subset I_a.\mathbf{P}_{\underline{i}}$$
Ce que nous souhaitions.\\

-Soit $n_a=n$. Puisque $I_a$ est $\mathfrak{m}_{\mathbf{R}}$-primaire et $\mathbf{R}/I$ n'est pas Artinien, on a  une inclusion stricte : $I\subset I_a$. Par cons\'equent il existe des $a_{i,j}\in \mathbf{R}$ tels que :
$$f_i=\sum_{j=1}^na_{i,j}g_j .$$
On a n\'ecessairement $det(a_{i,j})\in \mathfrak{m}_{\mathbf{R}}$, sans quoi $I=I_a$. Maintenant on a :
$$J_X(f)=det(a_{i,j}).J_X(g).$$
On conclut alors par 3.1 puisque $g_1,\ldots,g_n$ est r\'eguli\`ere et $det(a_{i,j})\in \sqrt{(g_1,\ldots,g_n)}=\mathfrak{m}_{\mathbf{R}}$.\endproof

\begin{ex}\textrm{ }\\
1) Si $\mathbf{A}=k$ un corps de caract\'eristique z\'ero, on retrouve l'implication $1)\Longrightarrow 2)$  du th\'eor\`eme 1.1.\\
2) Soit toujours $k$ un corps de caract\'eristique z\'ero, on peut prendre pour $\mathbf{A}=k[[Y_1,\ldots,Y_p]]/(h_1,\ldots,h_p)$ o\`u $(h_1,\ldots,h_p)$ engendre un id\'eal $(Y)$-primaire. Si $k$ est  muni d'une valeur absolue non triviale, on peut prendre de la m\^eme fa\c{c}on $\mathbf{A}=\mathcal{O}_p/(h_1,\ldots,h_p)$ o\`u $\mathcal{O}_p$ d\'esigne l'anneau des germes de s\'eries covergentes en $p$ variables et $(h_1,\ldots,h_p)$ est un id\'eal primaire pour l'id\'eal maximal de $\mathcal{O}_p$.\\
3) Consid\'erons toujours $S=k[[Y_1,\ldots,Y_p]]$ ou sa version convergente $\mathcal{O}_p$ dans le cas o\`u $k$ est muni d'une valeur absolue non triviale. Soit $\Delta $ un op\'erateur diff\'erentiel \`a coefficients constants :  
$$\Delta =\sum_{\alpha \in \mathbb{N}^p,\hspace{1pt} \vert \alpha \vert \leq m}a_{\alpha} \frac{\partial ^{\alpha }}{\partial Y^{\alpha }}$$
Alors $\Delta $ d\'efini un morphisme $\Delta _0 \in Hom_k(S,k)$ en posant $\Delta _0(f)=\Delta (f)(0)$.
Soit $J$ l'annulateur de $\Delta _0$, i.e. 
$$J=\{a \in S/ \forall f \in S,\hspace{3pt} \Delta _0(a.f)=0\}$$
Alors $J$ est un id\'eal primaire et $\mathbf{A}=S/J$ est un anneau de Gorenstein Artinien d'\'egale caract\'eristique z\'ero et le r\'esultat s'applique. Prenant par exemple $\Delta$  le \guillemotleft Laplacien \guillemotright : $\Delta =\sum_{i=1}^p\frac{\partial ^2}{\partial Y_{i}^{2}}$, on obtient un anneau de Gorenstein qui n'est pas intersection compl\`ete et qui n'est pas du type de l'exemple 2).
\end{ex}
\begin{cor}\textrm{ }\\
Soient $\mathbf{A}$ un anneau local de Gorenstein de dimension $0$ d'\'egale caract\'eristique z\'ero, $\mathbf{R}=\mathbf{A}[[X_1,\ldots,X_n]]$ et $f_1,\ldots,f_n$ des \'el\'ements de $\mathbf{R}$ tels que l'id\'eal $I$ engendr\'e par $f_1,\ldots,f_n$ soit un id\'eal propre. Alors les propri\'et\'es suivantes sont \'equivalentes :\\
1) $Dim \mathbf{R}/I>0$\\
2) $J_X(f_1,\ldots,f_n)\in I$
\end{cor}
\noindent \textit{Preuve :}
On vient de prouver $1)\Longrightarrow 2)$. $2)\Longrightarrow 1)$ ou de mani\`ere \'equivalente $non 1)\Longrightarrow non 2)$ d\'ecoule de $3.1)$. Puisque $R$ est de dimension $n$ et de Cohen-Macauley, si $(f_1,\ldots,f_n)$ ne satisfont pas $1)$ alors $I$ est primaire pour l'id\'eal maximal de $\mathbf{R}$ et donc la suite $(f_1,\ldots,f_n)$ est r\'eguli\`ere. On est donc dans les conditions alternatives de 3.1), comme nous l'avons plus haut.\endproof
\begin{rem}
M\^eme dans le cas le plus simple, c'est \`a dire $\mathbf{A}=\mathbb{C}$ et $\mathbf{R}=\mathcal{O}_n$ ou $\mathbb{C}[[X]]$ le r\'esultat n'est pas imm\'ediat \`a prouver par des techniques purement analytiques. Nous renvoyons \`a [B-Y] pour de telles tentatives.
\end{rem}

\section{Preuve de 1.2 et corollaire}
Par commodit\'e pour le lecteur, nous rappelons d'abord ici quelques notions essentielles pour lesquelles nos sources sont [H-I-O],[H-S], [L-T], [N-R], [R] et nous renvoyons \`a ces r\'ef\'erences pour de plus amples pr\'ecisions. Nous donnons ensuite une preuve du crit\`ere m\'etrique de d\'ependance int\'egrale de [L-T] qui n'utilise pas la d\'esingularisation et permet son extension en caract\'eristique arbitraire ou au contexte de tout corps muni d'une valeur absolue non triviale. De m\^eme, nous prouvons que dans le contexte des anneaux locaux d'\'egale caract\'eristique le crit\`ere valuatif de d\'ependance int\'egrale peut \^etre v\'erifier avec des arcs.  Nous d\'emontrons ensuite le th\'eor\`eme 1.2
\begin{propd}
Soient $\mathbf{R}$ un anneau commutatif untaire et noeth\'erien, $I$ un id\'eal de $\mathbf{R}$.
\begin{itemize}
\item[1)] Un \'el\'ement $x \in \mathbf{R}$ est dit entier sur $I$ si et seulement si il satisfait une relation :
$$x^k+a_1x^{k-1}+\ldots+a_k=0,\hspace{3pt}avec \hspace{3pt}a_i \in I^i,\hspace{3pt}1\leq i \leq k.$$
L'ensemble de tous les \'el\'ements entiers sur $I$ est un id\'eal $\overline{I}$ de $\mathbf{R}$, appel\'e cl\^oture int\'egrale de $I$.
\item[2)] Soit  $\theta =\frac{p}{q}$ un nombre rationnel positif, on dit qu'un \'el\'ement $x\in \mathbf{R}$ est dans $\overline{I^{\theta }}$ si et seulement si $x$ satisfait une relation :
$$x^k+a_1x^{k-1}+\ldots+a_k=0,\hspace{3pt}avec \hspace{3pt}\frac{ord_I(a_i)}{i}\geq \theta ,\hspace{3pt}1\leq i \leq k.$$
o\`u $ord_I(f)=Sup\{k\in \mathbb{N}\cup \{+\infty \}/f \in I^k\}$. Il revient au m\^eme de dire que $f^q\in \overline{I^p}$. $\overline{I^{\theta }}$ est un id\'eal de $\mathbf{R}$.
\end{itemize}
\end{propd}
\begin{lem}\textrm{ }\\
Soient $\mathbf{R}$ et $I$ comme ci-dessus.
\begin{itemize}
\item[1)] Pour toute partie multiplicativement ferm\'ee $S\subset \mathbf{R}$, $\overline{I}.\mathbf{R}_{S}=\overline{I.\mathbf{R}_S}$.
\item[2)] Si $\mathbf{R}\longrightarrow \mathbf{R}'$ est un morphisme fid\`element plat d'anneaux noeth\'eriens alors :
$$\overline{I.\mathbf{R}'}\cap \mathbf{R}=\overline{I}.$$
\end{itemize}
\end{lem}
L'appartenance d'un \'el\'ement $f$ \`a $\overline{I}$ se v\'erifie souvent \`a l'aide du crit\`ere valuatif de d\'ependance int\'egrale.
\begin{theo}(Crit\`ere valuatif de d\'ependance int\'egrale)([H-I-O],[H-R], [R])\\
Soient $\mathbf{R}$ un anneau noeth\'erien, $I\subset \mathbf{R}$ un id\'eal, $f \in \mathbf{R}$. Les propri\'et\'es suivantes sont \'equivalentes :\\
1) $f\in \overline{I}$\\
2) Pour toute valuation discr\`ete $v$: $\mathbf{R}\longrightarrow \mathbb{N}\cup \{ {+\infty}\}$, on a :
$$ v(f)\geq v(I),\hspace{3pt}v(I)=Min(v(x),x \in I)$$
De plus si $I$ n'est contenu dans aucun des premiers minimaux de $R$, il existe un ensemble fini de valuation discr\`ete $v_{l}$, $l \in \Lambda $, de rang $1$, unique \`a l'\'equivalence pr\`es des valuations et minimal pour la propri\'et\'e  : $f \in \overline{I}$ si et seulement si $v_{l}(f)\geq v_{l}(I)$, $\forall l \in\Lambda $. Cet ensemble de valuation est appel\'e l'ensemble des valuations de Rees de $I$.  Pour tout $k\in \mathbb{N}^*$, l'ensemble des valuations de Rees de $I$ et de $I^k$ coincident (c.f. [H-S] chap. 10). De plus pour tout $\theta \in \mathbb{Q}_{>0}$, $f\in \mathbf{R}$, on a $f \in \overline{I^{\theta }}$ si et seulement si on a :
$$v_l(f)\geq \theta .v_l(I),\hspace{3pt}\forall l \in \Lambda $$
\end{theo}
Les valuations de Rees peuvent \^etre construites par le proc\'ed\'e suivant.(c.f. [H-S] chap 10 Ex. 10.6 p.208). C'est une  \guillemotleft  version faisceautisée \guillemotright  de cette construction qui est adopt\'ee dans [L-T] dans le cadre de la g\'eom\'etrie analytique complexe. \\
Notons $\Pi '$ : $X'\longrightarrow Spec(\mathbf{R})$ l'\'eclatement normalis\'e de centre $I$. C'est \`a dire $\Pi '$ est le compos\'e de l'\'eclatement de centre $I$, $\Pi $ : $X\longrightarrow Spec(\mathbf{R})$, $X=Proj(\oplus _{m\geq 0}I^m)$, et du morphisme de normalisation $N$ : $X'\longrightarrow X$. Posons $Z=Spec(\mathbf{R}/I)$ et soit $Y'=\Pi '^{-1}(Z)$.  Consid\'erons la d\'ecomposition en composantes irr\'eductibles de $Y'_{red}$ i.e. :
$$Y'_{red}=\cup _{l \in \Lambda }Y'_l$$
Puisque $X'$ est normal, il est r\'egulier en codimension 1. Chaque $Y'_l$ \'etant irr\'eductible de codimension 1, notant $\mathcal{O}_{X',Y'_l}$ la fibre du faisceau structural  $\mathcal{O}_{X'}$ de $X'$ au point correspondant \`a l'id\'eal premier d\'efinissant $Y'_l$, chaque anneau $\mathcal{O}_{X',Y'_l}$ est un anneau r\'egulier de dimension 1 donc un anneau de valuation discr\`ete dont on note $v_l$ la valuation. On a pour tout $a\in \mathcal{O}_{X',Y'_l}$ :
$$v_l(a)=long(\mathcal{O}_{X',Y'_l}/a.\mathcal{O}_{X',Y'_l})$$
o\`u $long$ d\'esigne la longueur. Ou bien encore $v_l(a)=Max\{k\in \mathbb{N}/a \in \mathfrak{m}_{X',Y'_l}^k\}$
, o\`u l'on a d\'esign\'e par $\mathfrak{m}_{X',Y'_l}$ l'id\'eal maximal (principal) de $\mathcal{O}_{X',Y'_l}$. On fait alors agir les $v_l$ sur $R$ via le morphisme $\Pi '^{\sharp}$: $\mathbf{R}\longrightarrow \mathcal{O}_{X'}\longrightarrow \mathcal{O}_{X',Y'_l}$. On a ainsi une construction des valuations de Rees. Pour plus d'informations et une pr\'esentation dans le vocabulaire de l'alg\`ebre de Rees nous renvoyons \`a [H-S] chap. 10.
\begin{propd}(fonction asymptotique de Samuel) ([H-S],[S],[L-T])\\
Soient $\mathbf{R}$ et $I$ comme pr\'ec\'edemment, on appelle fonction asymptotique de Samuel de $I$ et on note $\overline{v}_I$ la fonction  : 
$$ \overline{v}_{I}(x)=Lim_{m\rightarrow +\infty}\frac{ord_I(x^m)}{m}=Sup_{\{m \in \mathbb{N}^*\}}\frac{ord_I(x^m)}{m}$$
qui est \`a valeurs dans $\mathbb{Q}^+\cup \{+\infty\}$.
On a : 
$$\overline{v}_{I}(x)= Min\{ \frac{v_{l}(x)}{v_{l}(I)}, \textit{ o\`u }v_{l} \textit{ parcourt l'ensemble des valuations de Rees de I} \}$$
Ainsi : $x \in \overline{I}$, si et seulement si $\overline{v}_{I}(x)\geq 1$. De fa\c{c}on plus g\'en\'erale, si $\theta \in \mathbb{Q}_{>0}$, on a : $x \in \overline{I^{\theta }}$ si et seulement si $\overline{v}_I(x)\geq \theta $.
\end{propd}
Dans [L-T], M. Lejeune et B. Teissier ont \'etudi\'e la notion de cl\^oture int\'egrale des id\'eaux en relation avec certaines in\'egalit\'es de \L ojasiewicz dans le cadre de la g\'eom\'etrie analytique complexe. Leurs preuves utilisent l'existence de d\'esingularisa\-tion et des propri\'et\'es topologiques de $\mathbb{C}$. Au moins deux faits importants ressortent de [L-T]. Le premier est le crit\`ere m\'etrique de d\'ependance int\'egrale qui permet de relier la fonction asymptotique de Samuel et certains exposants de \L ojasiewicz, le second est que l'on puisse v\'erifier le crit\`ere valuatif uniquement \`a l'aide d'arcs. Nous indiquons ici comment leurs r\'esultats se g\'en\'eralisent en caract\'eristique arbitraire et pour des corps alg\'ebriquement clos arbitraires (notamment en g\'eom\'etrie analytique rigide ) et comment la th\'eorie g\'en\'erale permet d'en donner une preuve.\\

L' essence du crit\`ere m\'etrique de d\'ependance int\'egrale est le lemme suivant.
\begin{lem}\textrm{ }\\
Soit $A$ un anneau (n\'ecessairement int\`egre) muni d'une valeur absolue non triviale : $\vert\hspace{3pt} \vert$, $A\longrightarrow [0,\infty[$. Alors si $a \in A$ satisfait une \'equation :
$$a^m+\sum_{i=1}^ma_ia^{m-i}=0$$
Alors : $\vert a\vert \leq 2.Max \vert a_i \vert^{\frac{1}{i}}$
\end{lem}
\noindent \textit{Preuve :}
Soit $K=frac{A})$ et $\hat{K}$ son compl\'et\'e pour la topologie d\'efinie par la valeur absolue. Alors la valeur absolue $\vert \hspace{3pt} \vert$ se prolonge de mani\`ere unique \`a $\hat{K}$. De m\^eme d\'esignons par $K'$ une cl\^oture alg\'ebrique de $\hat{K}$, alors il existe une unique extension de la valeur absolue de $\hat{K}$ \`a $K'$ (c.f [B-G-R]). On peut ainsi supposer sans perte de g\'en\'eralit\'e que $A$ est un corps alg\'ebriquement clos muni d'une valeur absolue non triviale. Le lemme est alors \'el\'ementaire. Soit en effet $i_0$ tel que $\vert a_{i_0}\vert^{\frac{1}{i_0}}=Max \vert a_i \vert^{\frac{1}{i}}$
et $\alpha \in k$ tel que $\alpha ^{i_0}=a_{i_0}$. Alors $a'=\frac{a}{\alpha }$ satisfait une \'equation :
$$(a')^m+\sum_{i=1}^mb_i (a')^{m-i}, \textrm{ avec }Max_{1\leq i \leq m} \vert b_i \vert^{\frac{1}{i}}\leq 1$$
Posons $r=\vert a'\vert$. Alors : $1\leq \sum_{i=1}^m \frac{1}{r^i}<\sum_{i=1}^{+\infty}\frac{1}{r^i}$
Donc n\'ecessairement $\frac{1}{r}\geq \frac{1}{2}$ et donc $r\leq 2$.\endproof

\begin{theo}\textrm{ }\\
Soit $k$ un corps alg\'ebriquement clos, $\mathbf{R}=k[[X_1,\ldots,X_n]]$ ou bien $\mathbf{R}=\mathcal{O}_n$ l'anneau des germes de s\'eries convergentes en $n$ variables, lorsque $k$ est munie d'une valeur absolue non triviale $\vert \quad \vert$ : $k\longrightarrow [0,+\infty[$. Soit $I=(u_1,u_2,\ldots,u_m)$ un id\'eal de $\mathbf{R}$. Alors pour $u\in \mathbf{R}$  et tout $\theta \in \mathbb{Q}_{>0}$ les propri\'et\'es suivantes sont \'equivalentes.\\
1) $u \in \overline{I^{\theta }}$.\\
2) (Crit\`ere m\'etrique de d\'ependance int\'egrale) Si $k$ est munie d'une valeur absolue non triviale et $\mathbf{R}=\mathcal{O}_n$ : 
$$\exists \varepsilon >0 ,\exists c>0 / \quad  x\in k^n et \parallel  x \parallel <\varepsilon \Longrightarrow (\vert u(x)\vert\leq C. (Max_{1\leq i \leq m}\vert u_i(x)\vert)^\theta $$ 
3) (k alg\`ebriquement clos quelconque, Crit\`ere valuatif par les arcs) Pour tout arc $\varphi ^*:$ $\mathbf{R}\longrightarrow k[[t]]$ 
$$ord_t(\varphi ^*(u))\geq \theta .Min( ord_t( \varphi^*(u_i)))$$
Il existe de plus un nombre fini d'arcs $\varphi ^{l*}$, $l \in \Lambda $ tels que $3)$ est satisfaite si et seulement si elle satisfaite pour ces arcs  $\varphi ^{l*}$, $l \in \Lambda $. On a de plus :
$$\overline{v}_I(u)=Min_{l\in \Lambda }\frac{ord_t(\varphi ^{l*}(u)}{Min_{1\leq i \leq n}\{ord_t(\varphi ^{l*}(u_i))\}}$$
\end{theo}
La donn\'ee d'un arc $\varphi ^*$: $\mathbb{R}\longrightarrow k[[t]]$, c'est \`a dire d'un morphisme n\'ecessairement local de $k$-alg\`ebres de $\mathbf{R}$ dans $k[[t]]$ est \'equivalente \`a la donn\'ee d'un \'el\'ement $\varphi (t) \in (tk[[t]])\times \ldots \times(tk[[t]])$ qui agit par composition sur $\mathbf{R}$ et nous emploierons indiff\'erement les deux notations. Le fait que l'on puisse v\'erifier le crit\`ere valuatif \`a l'aide d'arcs est une sp\'ecifit\'e due au fait que l'on travaille avec des anneaux locaux noeth\'eriens qui contiennent des corps (contexte d'\'egale caract\'eristique), contexte qui est suffisant pour nous. On peut faire des \'enonc\'es similaires en g\'eom\'etrie analytique rigide ou bien consid\'erer des alg\`ebres essentiellement de type fini sur un corps alg\'ebriquement clos.\\
\noindent \textit{Preuve : }(Sch\'ema)\\
$1)\Longrightarrow 2)$ par le lemme ci-dessus. $1)\Longrightarrow 3)$ imm\'ediat. $2)\Longrightarrow 3)$ car le corps $k$ contient des \'el\'ements arbitrairement petits pour la valeur absolue. $3)\Longrightarrow 1)$ r\'esulte facilement de la th\'eorie g\'en\'erale et du th\'eor\`eme de structure (I.S. Cohen) des anneaux complets r\'eguliers d'\'egales caract\'eristiques (En imitant [L-T]). En effet consid\'erons toujours $\Pi '$ : $X'\underrightarrow{N} X \underrightarrow{\Pi } Spec(R)$ l' \'eclatement normalis\'e de centre $I$ et $Z\hookrightarrow Spec(R)$, $Z=Spec(R/I)$ et toujours $Y'=\Pi'^{-1} (Z)$ et sa d\'ecomposition en composantes irr\'eductibles $Y'_{red}=\cup _{l \in \Lambda }Y'_{l}$ et soit : 
$$m'_{l}=long(\mathcal{O}_{Y',Y'_{l}})=long(\mathcal{O}_{X',Y'_{l}}/I.\mathcal{O}_{X',Y'_{l}})$$
$X'$ \'etant normal et donc r\'egulier en codimension 1, chaque $\mathcal{O}_{X',Y'_{l}}$ est un anneau r\'egulier de dimension 1 et donc un anneau de valuation discr\`ete dont on note $v_{l}$ la valuation. Ces valuations sont une description des valuations  de Rees comme nous l'avons mentionn\'e plus haut.  Posons $v_{l}(u.\mathcal{O}_{X',Y'_{l}})=v_{l}$. Notons $\mathcal{P}_{l}$ le sous-faisceau d'id\'eaux premiers  d\'efinissants $Y'_{l}$. D\'esignant par $\hat{ }$ la compl\'etion pour l'id\'eal maximal, le th\'eor\`eme de structure des anneaux r\'eguliers complets locaux d' \'egales caract\'eristiques permet d'affirmer que si $y'_{l}$ est un point ferm\'e, $y'_{l}\in Reg(Y'_{l})\cap  Reg(X')$ suffisament g\'en\'eral, on a en d\'esignant par $K(Y_l)$ le corps r\'esiduel de $\mathcal{O}_{X',Y_l}$:
\begin{align}
\hat{\mathcal{O}}_{X',Y'_{l}} &  \simeq K(Y'_{l})[[t]] \\
I.\hat{\mathcal{O}}_{X',Y'_{l}} &  \simeq  (t)^{m'_{l}}K(Y'_{l})[[t]]\\
\hat{\mathcal{O}}_{X',y'_{l}}  & \simeq  k[[t,U]], \quad U=(U_2,\ldots,U_n)\\
I. \hat{\mathcal{O}}_{X',y'_{l}} & \simeq  (t)^{m'_{l}}k[[t, U]]\\
u. \hat{\mathcal{O}}_{X',y'_{l}}   & \simeq  (t)^{v_{l}}k[[t,U]] \\
\mathcal{P}_{i',y'_{l}}.\hat{\mathcal{O}}_{X',y'_{l}} & \simeq (t).k[[t,U]]
\end{align}
L'arc :
 $$\varphi ^{l*}: \hspace{3pt}R \overset{\Pi'^{\sharp}}{\longrightarrow} \hat{\mathcal{O}}_{X',y'_{l}} =k[[t,U]]\longrightarrow k[[t]]$$
o\`u la derni\`ere fl\`eche est $t\rightarrow t$ et $U\longrightarrow 0$, calcule alors $v_{l}(I)$ (resp. $v_{l}(u)$) comme $Min_{ord_t}(\varphi ^{l*}(u_j))$ (resp. $ord_t(\varphi ^{l*})(u)))$) . Ainsi si 3) est satisfaite la consid\'eration des arcs $\varphi^{l*}$, $l \in \Lambda  $, nous assure que $v_{l}(u)\geq \theta .v_{l}(I)$, pour les valuations de Rees et donc que $u \in \overline{I^{\theta }}$ au vu des r\'esultats g\'en\'eraux rappel\'es ci-dessus. Il en va de m\^eme pour le calcul de la fonction asymptotique de Samuel. \endproof \\

\noindent \textit{Preuve de 1.2 :}\\
Soit $m$ le rang de la matrice Jacobienne $\frac{\partial f_i}{\partial X_j}(0)$. Alors le th\'eor\`eme des fonctions implicites formel nous assure qu'il existe un $k$-automorphisme $\alpha ^*$ de $k[[X_1,\ldots,X_n]]$ tel que :
$$\alpha ^*(f_i)=X_i,\hspace{3pt}1\leq i\leq m, \hspace{3pt}\alpha ^*(f_{m+i})=F_i,\hspace{3pt}1\leq i\leq n-m,$$
avec posant $X'=(X_{m+1},\ldots,X_n)$, $d_{X'}F_i(0)=0$, $1\leq i \leq n-m$. Posons $g_i(X')=F_i(0,X')$, $1\leq i \leq n-m$. On a $J_X(\alpha ^*(f_1),\ldots,\alpha ^*(f_n))=J_{X'}(g_1,\ldots,g_{n-m})$. Posons : 
\begin{align}I & = & (X_1,\ldots,X_m,F_1,\ldots,F_{n-m})& = &(X_1,\ldots,X_m,g_1(X'),\ldots,g_{n-m}(X'))\\
 J & = &(g_1,\ldots,g_{n-m})\subset k[[X']]
\end{align}
Pour $\theta \geq 1$, on a $J_X(f)\in \overline{I^{\theta }}$ si et seulement si $J_{X'}(g)\in \overline{J^\theta }$. Ainsi quitte \`a substituer les $g_i$ aux $f_i$, on pourra supposer sans perte de g\'en\'eralit\'es que $f_1,\ldots,f_n$ sont $n$ \'el\'ements de $k[[X_1,\ldots,X_n]]$ tels que le rang de la matrice Jacobienne $\frac{\partial f_i}{\partial X_j}(0)$ est nul. Il s'agit de voir que si $n=2$, alors $J_X(f)\in \overline{I}$ et que si $n>2$, il existe $\theta >1$ tel que $J_X(f)\in \overline{I^{\theta }}$. \\
Pour cela on peut supposer que $k$ est alg\'ebriquement clos. En effet si cela n'est pas le cas, soit $\overline{k}$ une cl\^oture alg\'ebrique de $k$. Le morphisme naturel $k[[X_1,\ldots,X_n]]\longrightarrow \overline{k}[[X_1,\ldots,X_n]]$ est fid\`element plat. Ainsi si $h\in k[[X]]$ et $I$ est un id\'eal on a : $h\in \overline{I^{\theta}}$ si et seulement si $h\in \overline{(I.\overline{k}[[X]])^{\theta }}$. D'autre part si $k$ est munie d'une valeur absolue non triviale $\vert \quad \vert$, on compl\`ete celui-ci pour la valeur absolue en un corps $k'$ et celle-ci se prolonge de mani\`ere unique \`a une cl\^oture alg\'ebrique de $\overline{k'}$  (c.f. [B.G.R]). Pour les questions ayant trait aux in\'egalit\'es, il suffit de v\'erifier celles-ci au voisinage de $0$ dans $\overline{k'}^n$. Nous supposerons donc que $k$ est alg\'ebriquement clos. D'apr\`es les r\'esultats mentionn\'es plus haut, notant $I=(f_1,\ldots,f_n)$, il existe un nombre fini d'arcs $\varphi ^{l}$, $l\in \Lambda $, $\varphi ^{l} \in (tk[[t]])\times\ldots \times (tk[[t]])$ tels que pour tout $\theta=p/q \in \mathbb{Q}_{>0}$ et tout $h \in k[[X]]$ on a : $h\in \overline{I^{\theta }}$ si et seulement si :
$$\forall l \in \Lambda , ord_t(h\circ \varphi ^l)\geq \theta. Min_{1\leq i \leq n} ord_t(f_i\circ \varphi ^l)$$
Fixons un tel arc $\varphi $ (nous omettrons l'exposant $l$). On a alors :
$$(f_i\circ \varphi )'(t)=\sum_{j=1}^n \frac{\partial f_i}{\partial X_j}(\varphi (t))\varphi _j'(t)$$
Par les r\`egles de Cramer, d\'esignant par $\Delta _{i,j}$ le mineur obtenu \`a partir de la matrice Jacobienne en rayant la $i^{\grave{e}me}$ ligne et la $j^{\grave{e}me}$ colonne on a :
$$\varphi' _j(t). J_X(f)(\varphi (t)=\sum_{i=1}^n (-1)^{i+j}\Delta _{i,j}(\varphi (t))(f_i\circ \varphi)'(t) $$
Comme la caract\'eristique de $k$ est z\'ero, on a $ord_t(\varphi '_j)=ord_t(\varphi _j)-1$ et $ord_t(f_i\circ \varphi )'=ord_t(f_i\circ \varphi ) -1$. Ainsi l'identit\'e ci-dessus fournit :
$$ord_t(\varphi _j)+ord_t(J_X(f)(\varphi )\geq Min_{1 \leq i \leq n}(ord_t(\Delta _{i,j}(\varphi (t))+ord_t(f_i\circ \varphi ))$$
Fixons alors $j$ tel que $ord_t (\varphi _j)=Min_{1\leq k\leq n}ord_t(\varphi _k)=ord_t(\varphi ^*(\mathfrak{m}))=ord_t(\varphi )$ o\`u $\mathfrak{m}$ est l'id\'eal maximal de $k[[X]]$. Pour un tel indice $j$, on obtient par l'in\'egalit\'e ci-dessus :
$$ord_t(J_X(f)(\varphi )\geq Min_{1\leq i\leq n}(ord_t(\Delta _{i,j}(\varphi ))-ord_t(\varphi ))+ Min_{1\leq i \leq n} (ord_t(f_i\circ \varphi ))$$
Maintenant puisque le rang de la matrice jacobienne de $(f_1,\ldots,f_n)$ en $0$ est nul, on a : $\Delta _{i,j}\in \mathfrak{m}$ si $n=2$ et $\Delta _{i,j} \in \mathfrak{m}^{2}$ si $n\geq 3$. Par suite l'in\'egalit\'e ci-dessus assure que $J_X(f)\in \overline{I}$ si $n=2$. Si $n\geq3$, alors puisque $\Delta _{i,j}\in \mathfrak{m}^2$ on a :
$$ ord_t(\Delta _{i,j}(\varphi ))-ord_t(\varphi )\geq ord_t(\varphi )$$
Maintenant $I$ est un id\'eal $\mathfrak{m}$-primaire, donc il existe $s\in \mathbb{N}^*$ tel que : $\mathfrak{m}^s \subset I$. D'o\`u $s.ord_t(\varphi )\geq  Min_{1\leq i \leq n} ord_t(f_i\circ \varphi )$.
Ce qui permet d'obtenir :
$$ord_t(J_X(f)(\varphi )\geq (1+1/s)Min_{1\leq i \leq n}ord_t(f_i\circ \varphi )$$
Ainsi le nombre rationnel $\theta $ :
$$\theta =Min_{l\in \Lambda }\frac{ord_t(J_X(f)\circ \varphi ^l)}{Min_{1\leq i\leq n}(ord_{t}(f_i\circ \varphi ^l))}=\overline{v_{I}}(J_X(f))$$
est strictement plus grand que $1$. Les in\'egalit\'es de \L ojasiewicz de 1.2) dans le cas convergent s'obtiennent par le crit\`ere m\'etrique de d\'ependance int\'egrale mentionn\'e plus haut.\endproof

\begin{rem}\textrm{ }\\
1) Si $n-m=2$, il n'existe pas en g\'en\'eral de $\theta >1$ plus grand que $1$ tel que : $J_X(f_1,\ldots,X_n)\in \overline{I^{\theta }}$. Il suffit de consid\'erer l'exemple :
$$(f_1,\ldots,f_n)=(X_1,\ldots,X_{n-2},X_{n-1}^2,X_n^2).$$
On a alors $J_X(f) \in \overline{I}$ mais pour aucun nombre $\theta >1$, $J_X(f_1,\ldots,f_n)\in \overline{I^{\theta }}$. En effet ici $\overline{I}=(X_1,\ldots,X_{n-2},X_{n-1}^2,X_{n-1}X_n,X_{n}^2)$ et $J_X(f)=4X_{n-1}X_{n}$. Mais c'est \`a peu pr\`es le seul cas o\`u une telle  situation se produit. La m\^eme preuve que ci-dessus montre que $\overline{v_{I}}(J_X(f))>1$ d\`es que $Min(ord_X(f_i))\geq 3$.\\
2) Bien entendu le cas de $k=\mathbb{R}$ ou $\mathbb{C}$ est le premier qui se pr\'esente. Mais on peut consid\'erer aussi le cas o\`u $R$ est un corps r\'eel clos et prendre pour $k$ le corps des s\'eries de puiseux \`a coefficients dans $R$ ou $C=R+iR$ i.e $k=\cup _{q \in \mathbb{N}^*} R((X^{1/q}))$ ou $k=\cup _{q \in \mathbb{N}^*}C((X^{1/q}))$. De tels corps sont naturellement munis d'une valeur absolue ultram\'etrique non triviale $\vert a \vert=e^{-ord_X(a)}$. Les in\'egalit\'es de \L ojasiewicz de 1.2 sont valables pour les s\'eries convergentes \`a coefficients dans de tels corps.\\
3) Un c\'el\`ebre r\'esultat de S. \L ojasiewicz affirme que si $f$ est un germe de s\'erie convergente en $n$ variables r\'eelles, $f(0)=0$, alors il existe un nombre rationnel $\theta >1$ tel que au voisinage de $0$ dans $\mathbb{R}^n$, on ait $\vert f(x)\vert\leq C (\sum_{1\leq i \leq n} \vert \frac{\partial f}{\partial X_i}(x)\vert)^{\theta }$.
Pour des preuves tr\`es diff\'erentes les unes des autres nous renvoyons respectivement \`a [B-M] prop. 6.8 p.35, [T2] compl\'ement 1) et $[Loj]$. Ce r\'esultat peut \^etre \'etendu au cas d'un corps de caract\'eristique z\'ero muni d'une valeur absolue non triviale en proc\'edant comme ci-dessus. Les in\'egalit\'es de 1.2) qui nous ont \'et\'e inspir\'ees par [T1] sont les analogues pour le Jacobien de ces in\'egalit\'es. 
\end{rem}

Signalons pour terminer le corollaire imm\'ediat.

\begin{cor}\textrm{ }\\ 
Soient $k$ un corps alg\'ebriquement clos, de caract\'eristique z\'ero, muni d'une valeur absolue non triviale $\vert\quad \vert$: $k\longrightarrow [0,\infty[$ et $\mathcal{O}_n$ l'anneau des germes de s\'eries convergentes en $n$ variables $(z_1,\ldots,z_n)$.\\
 1)Soit $f\ \in  \mathcal{O}_{n}$, $f(0)=0$. Alors f est \`a singularit\'e isol\'ee \`a l'origine si et seulement si son Hessien $H_{Z}(f)$ n'appartient pas \`a l'id\'eal jacobien de $f$.\\ 
 2) Soit $f\ \in  \ \mathcal{O}_{n}$, $f(0)=0$ \`a singularit\'e isol\'ee \`a l'origine. Alors $f$ est, \`a
changement de variables pr\`es, \'egale \`a une singularit\'e du type $z_{1}^{2}+z_{2}^{2}+\cdots+z_{n}^{k}$, $k\geq 2,$ si
et seulement si son id\'eal jacobien est int\'egralement clos.
\end{cor}
\noindent \textit{Preuve :}\\
 1) est la traduction directe au cas consid\'er\'e des r\'esultats de 3). 2) Les  singula\-rit\'es du type consid\'er\'e ont bien \'evidemment un id\'eal jacobien int\'egralement clos. R\'eciproquement, soit $\mathcal{H}_{Z}(f)$  la matrice Hessienne de $f$. Soit $m$ le rang de  $\mathcal{H}_{Z}(f)$  en  $0$. Si $n-m\geq 2$
alors $H_{Z}(f)$ est dans la cl\^oture int\'egrale de l'id\'eal  $(\partial  f/\partial  z_{1},\ldots,\partial  f/\partial z_{n})$ . Par cons\'equent si celui-ci \'etait int\'egralement clos on aurait $H_{Z}(f)\in   (\partial  f/\partial  z_{1},\ldots,\partial f/\partial  z_{n})$. Ceci  est absurde en vertu du 1) et de l'hypoth\`ese de singularit\'e isol\'ee. Donc la seule possibilit\'e pour que l'id\'eal jacobien de $f$ soit int\'egralement clos est que le rang  $m$ de  $\mathcal{H}_{Z}(f)$ en z\'ero satisfasse $n-m <2$. Par suite $m\geq n-1$. Ceci implique que $f$ est du type propos\'e, en vertu du lemme de Morse \`a param\`etres (c.f. [A-G-V] p.187).\endproof

\vspace{10pt}

\centerline{ \textit{BIBLIOGRAPHIE }}

\vspace{10pt}

\noindent [A-G-V] V.I. ARNOLD, S.M. GUSEIN-SAIDE, A.N. VARSCHENKO, \textit{Singularities of differentiable Maps, vol 1}, Monographs in Mathematics vol. 82, Birkh\"auser 1985.

\noindent [B1] N. BOURBAKI, \textit{Alg\`ebre Chapitre 1 \`a 3}, nouvelle \'edition Hermann 1970.

\noindent [B2] N. BOURBAKI, \textit{Alg\`ebre commutative chap.1 \`a 7}, Hermann 1961. 

\noindent [B3] N. BOURBAKI, \textit{Alg\`ebre commutative chap. 10}, Masson 1998.

\noindent [B-G-R] S. BOSCH, U. GERRITZEN, R. REMMERT, \textit{Archimedian Analysis, A systematic Approach to Rigid Analytic Geometry}, Grundlehren der mathematichen Wissenschaften 261, Springer 1984.

\noindent [B-H1] J.Y. BOYER, M. HICKEL, \textit{Une g\'en\'eralisation de la loi de Transformation pour les r\'esidus}, Bull. Soc. Math. France $\mathbf{125}$, (1997), 315-335.

\noindent [B-H2] J.Y. BOYER, M. HICKEL, \textit{Extension dans un cadre alg\'ebrique d'une Formule de Cauchy-Weil}, Manuscripta Mathematica $\mathbf{98}$, (1999), 195-223.

\noindent [B-M] E. BIERSTONE, P.D. MILMAN, \textit{Semianalytic and Subanalytic Sets}, Pub. Math. I.H.E.S. $\mathbf{67}$,(1988), 1-42. 

\noindent [BO] J. Y. BOYER, \textit{Comparaison de diff\'erentes approches des r\'esidus et applications}, Th\`ese de l'universit\'e Bordeaux 1, F\'evrier 1999.

\noindent [B-Y] C.A. BERENSTEIN, A. YGER, \textit{Analytic residue theory in the non-complete intersection case}, J. Reine angew. Math. $\mathbf{527}$ (2000), 203-235.

 \noindent[E] D. EISENBUD, \textit{Commutative algebra with a view toward algebraic geometry}, Graduate Texts in Mathematics $n^{°}$ 150, Springer-verlag.

\noindent [Hi] M. HICKEL, \textit{Une remarque \`a propos du Jacobien de n fonctions holomorphes \`a l'origine de $\mathbb{C}^n$ }, accept\'e pour publication \` a Ann. Pol. Math. 20 Mars 2000 , retir\'e de la publication par l'auteur en Juin 2000.

\noindent [H-I-O]M. HERMAN, S. IKEDA, U. ORBANZ, \textit{Equimultiplicity and blowing up, an algebraic study with an appendix by B. MOONEN}, Springer 1988.

 \noindent [H-K] R. HUEBL, E. KUNZ, \textit{Integration of differential forms on schemes}, J. Reine Angew. Math. $\mathbf{410}$ (1990), 53-83.

 \noindent [H-S] C. HUNEKE-I. SWANSON, \textit{Integral Closure of ideals, Rings, and Modules}, London Mathematical Society Lecture Note Series $n^{°}$ 336, Cambridge University Press 2006.

\noindent [L] J. LIPMAN, \textit{Residues and traces of differential forms via Hochschild homology}, Contemporary Mathematics, vol. 61, American Mathematical Society, Providence-Rhode-island.

\noindent [Lo] S. \L OJASIEWICZ, \textit{Ensembles semi-analytiques}, Pub. Math. I.H.E.S. 1964.

\noindent [L-T] M. LEJEUNE-B. TEISSIER, \textit{Cl\^oture int\'egrale des id\'eaux et \'equisingularit\'e}, S\'eminaire Lejeune-Teissier, Centre de Math\'ematiques de l'\'ecole polytechnique 1974, Publications Universit\'e Scientifique et M\'edicale de Grenoble.

\noindent [M] H. MATSUMURA, \textit{Commutative Ring theory}, Cambridge studies in mathematics 8, 1986.

\noindent [N] E. NETTO, \textit{Vorlesungen \"uber algebra, vol. 2}, Leipzig, Teubner 1900.

\noindent [P] A. PLOSKI, \textit{Question orale lors du colloque Effectivity problems algebraic and analytic methods}, University of Calabria, 22-28 Juin 98.

\noindent [N.R] D.G. NORTHCOTT-D. REES, \textit{Reduction of ideals in local rings}, Proc. Cambridge Phil. Soc., $\mathbf{50}$ 2 (1954), 145-158.

\noindent [R] D. REES, \textit{Lectures on the asymptotic theory of ideals}, London Mathematical lecture notes series 113, Cambridge University Press.

\noindent [S] S. SPODZIEJA, \textit{On some property of the jacobian of a homogeneous polynomial mapping}, Bull. Soc. Sci et Lettres de Lodz, Vol XXXIX, 5 (1989), 1-5.

\noindent [S-S] G. SCHEJA, U. ST\" ORCH, \textit{\"Uber spurfunktionen bei vollst\"andigen durchschnitten}, J. Reine angew. Math. $\mathbf {278/279}$ (1975), 174-190.

\noindent [T1] B. TEISSIER, \textit{Cycles \'evanescents, sections planes et conditions de Whitney}, Ast\'erisque 7 et 8, (1973), 285-362.

\noindent [T2] B. TEISSIER, \textit{Sept compl\'ements au s\'eminaire Lejeune-Teissier}, \`a para\^{\i}tre dans Ann. Fac. Sci. Toulouse.

\noindent [W1] W.V. VASCONCELOS, \textit{The top of a system of equations}, Bol. Soc. Mat. Mexicana (issue dedicated to Jos\'e Adem), $\mathbf{37}$ (1992), 549-556.
 
\noindent[W2] W.V. VASCONCELOS, \textit{Arithmetic of Blowup Algebra}, London Mathematical Lecture Note Series $n^{°}$ 195, Cambridge University Press, 1994.

\noindent [W3] W.V. VASCONCELOS, \textit{Computational Methods in Commutative Algebra and Algebraic Geometry}, Algorithms and computation in mathematics volume 2, Springer-Verlag 1997.\\

 Michel HICKEL,\\ 
Universit\'e Bordeaux 1, I.M.B.\\
Equipe d'Analyse et G\'eom\'etrie\\
et I.U.T. Bordeaux 1 d\'epartement Informatique\\ 
33405 Talence Cedex, France\\
email : hickel@math.u-bordeaux1.fr 
\end{document}